\newif\ifreport
\reporttrue
\ifreport
  \documentclass{article}\def\zibreport{1}
\else
  \documentclass{llncs}\def\zibreport{0}
\fi

\usepackage{authblk}
\makeatletter
\renewcommand\AB@authnote[1]{\rlap{\textsuperscript{\normalfont#1}}}

\makeatother

\makeatletter
\DeclareRobustCommand*{\escapeus}[1]{%
  \begingroup\@activeus\scantokens{#1\endinput}\endgroup}
\begingroup\lccode`\~=`\_\relax
   \lowercase{\endgroup\def\@activeus{\catcode`\_=\active \let~\_}}
\makeatother

\usepackage[T1]{fontenc}
\usepackage[utf8x]{inputenc}
\usepackage[english]{babel}

\usepackage{booktabs}
\usepackage{pdflscape}
\usepackage{tabularx}
\usepackage{longtable}

\usepackage{tikz}
\usetikzlibrary{decorations.pathmorphing}

\usepackage{pgfplots, pgfplotstable}

\usepackage[font=small]{caption}
\usepackage[textsize=small]{todonotes}

\newcommand{\thetitle}{Experiments with Conflict Analysis in Mixed Integer Programming}

\usepackage{amsmath}
\usepackage{amssymb}

\usepackage{xspace}
\usepackage{xcolor}

\usepackage{ifthen}
\usepackage{xifthen}

\usepackage[colorlinks=false, pdftex]{hyperref}

\newcommand{\ie}{i.e.,\xspace}
\newcommand{\eg}{e.g.,\xspace}

\newcommand{\scip}{\textsf{SCIP}\xspace}
\newcommand{\soplex}{\textsf{SoPlex}\xspace}

\newcommand{\bandb}{branch-and-bound\xspace}
\newcommand{\Bandb}{Branch-and-bound\xspace}

\newcommand{\intvars}{\mathcal{I}}
\newcommand{\allvars}{\mathcal{N}}

\newcommand{\optMIP}[1]{#1^{\star}}
\newcommand{\optLP}[1]{#1^{\star}_{LP}}

\newcommand{\maxact}[4]{\ifthenelse{\isempty{#4}}{\Delta_{\max}(#1,#2,#3)}{\Delta_{\max}^{#4}(#1,#2,#3)}}

\newcommand{\Z}{\mathbb{Z}\xspace}
\newcommand{\R}{\mathbb{R}\xspace}

\newcommand{\conflict}{\texttt{conflict}\xspace}
\newcommand{\combined}{\texttt{combined}\xspace}
\newcommand{\dualray}{\texttt{dualray}\xspace}
\newcommand{\combinedpool}{\texttt{combined+pool}\xspace}

\newcommand{\Miplib}{\textsc{Miplib}}
\newcommand{\Conflict}{\textsc{Conflict}\xspace}

\ifthenelse{\zibreport = 1}{
  \usepackage{zibtitlepage}

  \ZTPTitle{\thetitle}
  \ZTPAuthor{Jakob Witzig\hspace{1cm} Timo Berthold\hspace{1cm} Stefan Heinz}
  \ZTPPreprint
  \ZTPNumber{16-63}
  \ZTPMonth{November}
  \ZTPYear{2016}
}{}

\begin{document}

\title{\thetitle}

\ifthenelse{\zibreport = 0}{
    \author{
      Jakob Witzig\inst{1}
      \and Timo Berthold\inst{2}
      \and Stefan Heinz\inst{2}
    }
    \institute{
      Zuse Institute Berlin, Takustr.~7, 14195~Berlin, Germany\\ \email{witzig@zib.de}
      \and
      Fair Isaac Germany GmbH, Takustr.~7, 14195~Berlin, Germany\\ \email{\{timoberthold,stefanheinz\}@fico.com}
    }
}{
  \author[1]{Jakob~Witzig}
  \author[2]{Timo~Berthold}
  \author[2]{Stefan~Heinz}

  \affil[1]{Zuse Institute Berlin, Takustr.~7, 14195~Berlin, Germany \protect\\ \texttt{witzig@zib.de}\medskip}
  \affil[2]{Fair Isaac Germany GmbH, Takustr.~7, 14195~Berlin, Germany \protect\\ \texttt{\{timoberthold,stefanheinz\}@fico.com}}
}

\ifthenelse{\zibreport = 1}{\zibtitlepage}{}

\maketitle

\begin{abstract}
  The analysis of infeasible subproblems plays an import role in
  solving mixed integer programs (MIPs) and is implemented in most
  major MIP solvers.  There are two fundamentally different concepts
  to generate valid global constraints from infeasible subproblems.
  The first is to analyze the sequence of implications obtained by
  domain propagation that led to infeasibility. The result of the
  analysis is one or more sets of contradicting variable bounds from
  which so-called conflict constraints can be generated.  This concept
  has its origin in solving satisfiability problems and is similarly
  used in constraint programming. The second concept is to analyze infeasible
  linear programming (LP) relaxations. The dual LP solution provides a
  set of multipliers that can be used to generate a single new
  globally valid linear constraint.
  The main contribution of this short paper is an empirical
  evaluation of two ways to combine both approaches.
  Experiments are carried out on general MIP instances from standard
  public test sets such as \Miplib2010; the presented algorithms have been
  implemented within the non-commercial MIP solver SCIP.  Moreover, we
  present a pool-based approach to manage conflicts which addresses
  the way a MIP solver traverses the search tree better than aging
  strategies known from SAT solving.
\end{abstract}

\section{Introduction: MIP and Conflict Analysis}
\label{sec:introduction}

In this paper we consider \emph{mixed integer programs (MIPs)} of the form
\begin{align}
  \optMIP{c} = \min \{ c^tx\,|\,Ax \geq b,\ \ell \leq x \leq u,\ x \in \Z^{k} \times \R^{n-k} \}, \label{eq:mip}
\end{align}
with objective function $c \in \R^{n}$, constraint matrix $A \in \R^{m\times n}$,
constraint left-hand side $b \in \R^{m}$, and variable bounds $\ell,u \in \overline{\R}^{n}$,
where $\overline{\R} := \R \cup \{\pm\infty\}$.
Furthermore, let $\allvars = \{1,\ldots,n\}$ be the index set of all variables.
Let $\intvars \subseteq \allvars$ such that $x_{i} \in \Z$ for all
$i \in \intvars$, \ie the set of variables that need to be integral in every feasible solution.

When omitting the integrality requirements, we obtain the \emph{linear program (LP)}
\begin{align}
  \optLP{c} = \min \{ c^tx\,|\,Ax \geq b,\ \ell \leq x \leq u,\ x \in \R^{n} \}. \label{eq:lprelax}
\end{align}
The linear program~\eqref{eq:lprelax} is called \emph{LP relaxation} of~\eqref{eq:mip}.
The LP relaxation provides a lower bound on the optimal solution value
of the MIP~\eqref{eq:mip},
\ie $\optLP{c} \leq \optMIP{c}$.
In LP-based \bandb~\cite{Dakin1965,LandDoig1960}, the most commonly
used method to solve MIPs, the LP relaxation is used for bounding.
\Bandb is a divide-and-conquer method which splits
the search space sequentially into smaller subproblems that are
(hopefully) easier to solve.
During this procedure we may encounter infeasible subproblems.
Infeasibility can be detected by contradicting implications, \eg derived by domain propagation, or
by an infeasible LP relaxation.
Modern MIP solvers try to learn from infeasible subproblems, e.g., by \emph{conflict analysis}.
Conflict analysis for MIP has its origin in solving satisfiability problems (SAT) and goes back to~\cite{marques1999grasp}. 
Similar ideas are used in constraint programming, \eg see~\cite{ginsberg1993dynamic,jiang1994no,stallman1977forward}. 
First integrations of these techniques into MIP were independently
suggested by~\cite{DaveyBolandStuckey2002},~\cite{SandholmShields2006}
and~\cite{achterberg2007constraint}.
Further publications suggested to use conflict information for variable
selection in branching, to
tentatively generate conflicts before
branching~\cite{AchterbergBerthold2009,KilincNemhauserSavelsbergh2010},
and to analyze infeasibility detected in  primal heuristics~\cite{BertholdGleixner2013,BertholdHendel2015}.

Today, conflict analysis is widely established in solving MIPs.
The principal idea of conflict analysis, in MIP terminology, can be sketched as follows.

Given an infeasible node of the \bandb tree defined by the subproblem
\begin{align}
  \min \{ c^tx\,|\,Ax \geq b,\ \ell^{\prime} \leq x \leq u^{\prime},\ x \in \Z^{k} \times \R^{n-k} \} \label{eq:infsubproblem}
\end{align}
with local bounds $\ell \leq \ell^{\prime} \leq u^{\prime} \leq u$.
In LP-based \bandb, the infeasibility of a subproblem is typically  detected by an infeasible LP
relaxation (see next section) or by contradicting implications.

In the latter case, a \emph{conflict graph} gets constructed which
represents the logic of how the set of branching decisions led to the
detection of infeasibility.  More precisely, the conflict graph is a
directed acyclic graph in which the vertices represent bound changes
of variables and the arcs $(v,w)$ correspond to bound changes implied
by propagation, \ie the bound change corresponding to $w$ is based
(besides others) on the bound change represented by $v$.  In addition
to these inner vertices which represent the bound changes from domain
propagation, the graph features source vertices for the bound changes
that correspond to branching decisions and an artificial sink vertex
representing the infeasibility.  Then, each cut that separates the
branching decisions from the artificial infeasibility vertex gives
rise to a valid \emph{conflict constraint}. A conflict constraint
consists of a set of variables with associated bounds, requiring that in each feasible solution
at least one of the variables has to take a value outside these bounds.
Note that in general, this is not a linear constraint and that by
using different cuts in the graph, several different conflict
constraints might be derived from a single infeasibility.
A variant of conflict analysis close to the one described above is
implemented in \scip, the solver in which we will conduct our computational experiments.
Also, a similar implementation is available in the \textsf{FICO Xpress-Optimizer}.

This short paper consists of two parts which are independent but complement each other in practice.
The first part of this paper (Section~\ref{sec:analyzingdualsols}) focuses on a MIP technique to analyze infeasibility based on LP theory.
We discuss the interaction, differences, and commonalities between
conflict analysis and the so-called \emph{dual ray analysis}.
Although both techniques have been known before, \eg~\cite{achterberg2007constraint,polik15ics}, 
this will be, to the best of our knowledge, the first published direct
comparison of the two.
In the second part (Section~\ref{sec:conflictpool}), we present a new approach to drop conflicts that do not lead to variable bound reductions frequently.
This new concept is an alternative to the aging scheme known from SAT.
Finally, we present computational experiments comparing the techniques described in Section~\ref{sec:analyzingdualsols} and~\ref{sec:conflictpool}.\\

\section{Analyzing Dual Unbounded Solutions}
\label{sec:analyzingdualsols}

The idea of conflict analysis is tightly linked to domain
propagation: conflict analysis studies a sequence of variable bound implications made by
domain propagation routines. Besides domain
propagation, there is another important subroutine in MIP solving
which might prove infeasibility of a subproblem: the LP relaxation.
The proof of LP infeasibility  comes in form of a so-called ``dual
ray'', that is a list of multipliers on the model constraints and the
variable bounds. Those give rise to a globally valid constraint that
can be used similarly to a conflict constraint. In this section, we
discuss the analysis of the LP infeasibility proof in more detail.

\subsection{Analysis of Infeasible LPs: Theoretical Background}
\label{subsec:theory}

Consider a node of the \bandb tree and the corresponding subproblem
\begin{align}
  \min \{ c^tx\,|\,Ax \geq b,\ \ell^{\prime} \leq x \leq u^{\prime},\ x \in \Z^{k} \times \R^{n-k} \} \label{eq:subproblem}
\end{align}
defined by local bounds $\ell \leq \ell^{\prime} \leq u^{\prime} \leq u$.
The \emph{dual LP} of the corresponding LP relaxation of~\eqref{eq:subproblem} is given by
\begin{align}
  \max \{ y^tb + \underline{r}^t\ell^\prime + \overline{r}^tu^\prime\,|\,A^ty + \underline{r} + \overline{r} \leq c,\ y,\, \underline{r} \in \R^{n}_{\geq 0},\, \overline{r} \in \R^{n}_{\leq 0} \}, \label{eq:duallp}
\end{align}
where $A_{\cdot i}$ is the $i$-th column of $A$, $\underline{r}_i = \max\{0, c_i - y^tA_{\cdot i}\}$, and $\overline{r}_i = \min\{0, c_i - y^tA_{\cdot i}\}$.
By LP theory each unbounded ray $(\gamma,\underline{r},\overline{r})$ of~\eqref{eq:duallp} proves infeasibility of~\eqref{eq:subproblem}.
A ray is called unbounded if multiplying the ray with an arbitrary scalar $\alpha > 0$ will not change the feasibility.
Note, in this case it holds
\begin{align*}
  \underline{r}_i = \max\{0, - y^tA_{\cdot i}\} \quad \text{and} \quad \overline{r}_i = \min\{0, - y^tA_{\cdot i}\}.
\end{align*}
Moreover, the Lemma of Farkas states that exactly one of the following two systems is satisfiable
\begin{align*}
  \left.
  (F_{1})\quad
  \begin{array}{rl}
    Ax&\geq b \\
    \ell^{\prime} \leq x& \leq u^{\prime}
  \end{array}
  \right\} \dot{\vee}
  \left\{
  \begin{array}{rl}
    \gamma^tA + \underline{r} + \overline{r}&\leq 0 \\
    \gamma^tb + \underline{r}^t\ell^{\prime} + \overline{r}^tu^{\prime} &> 0
  \end{array}
  \quad(F_{2})
  \right.
\end{align*}
It follows immediately, that if $F_{1}$ is infeasible, there exists an unbounded ray $(\gamma,\underline{r},\overline{r})$ of~\eqref{eq:duallp} satisfying $F_{2}$.
An infeasibility proof of~\eqref{eq:subproblem} is given by a single constraint
\begin{align}
  \gamma^tAx \geq \gamma^tb, \label{eq:farkascons}
\end{align}
which is an aggregation of all rows $A_{j\cdot}$ for $j = 1,\ldots,m$ with weight $\gamma_{j} > 0$.
Constraint~\eqref{eq:farkascons} is globally valid but violated in the local bounds $[\ell^{\prime},u^{\prime}]$ of subproblem~\eqref{eq:subproblem}.
In the following, this constraint will be called \emph{proof-constraint}.

\subsection{Conflict Analysis of Infeasible LPs}
\label{subsec:confanalysisinflp}

The analysis of an infeasible LP relaxation, as it is implemented in
\scip, is a hybrid of the theoretical considerations made in Section~\ref{subsec:theory} and the analysis of the conflict graph known from SAT.
To use the concept of a conflict graph, all variables with a non-zero
coefficient in the proof-constraint are converted to vertices of the conflict graph representing bound changes; global bound changes are omitted.
Those vertices, called the \emph{initial reason}, are then connected to the artificial sink representing the infeasibility.
This neat idea was introduced in~\cite{achterberg2007conflict}. From
thereon, conflict analysis can be applied as described in Section~\ref{sec:introduction}.

In practice, the proof-constraint is often quite dense, and therefore,
it might be worthwhile to search for a sparser infeasibility proof.
This can be done by a heuristic that relaxes some of the local bounds $[\ell^\prime,u^\prime]$ that appear in the proof-constraint.
Of course, the relaxed local bounds $[\ell^{\prime\prime},u^{\prime\prime}]$ with 
$\ell < \ell^{\prime\prime} \leq \ell^{\prime} \leq u^{\prime} \leq u^{\prime\prime} < u$ still need to fulfill
\begin{align*}
  \gamma^tb + \underline{r}^t\ell^{\prime\prime} + \overline{r}^tu^{\prime\prime} > 0.
\end{align*}
The more bounds can be relaxed that way, the smaller gets the initial
reason and consequently the stronger are the derived conflict constraints.
Note again that these constraints do not need to be linear, if general integer or continuous variables are present.

\subsection{Dual Ray Analysis of Infeasible LPs}
\label{subsec:dualrays}

The proof-constraint is globally valid but infeasible within the local bounds.
It follows immediately by the Lemma of Farkas that the \emph{maximal activity}
\begin{align*}
  \maxact{\gamma^tA}{\ell^{\prime}}{u^{\prime}}{} := \sum_{i \in \allvars\colon \gamma^tA_{i} > 0} (\gamma^tA_{i}) u^{\prime}_{i} 
  + \sum_{i \in \allvars\colon \gamma^tA_{i} < 0} (\gamma^tA_{i}) \ell^{\prime}_{i}
\end{align*}
of $\gamma^tAx$ w.r.t.\ variable bounds $[\ell^\prime,u^\prime]$ is strictly less than the corresponding left-hand side $\gamma^tb$.

 Instead of creating an ``artificial'' initial reason, the
 proof-constraint might also be used directly for domain propagation
 in the remainder of the search. It is a conic combination of global
 constraints, \ie it is itself a valid (but redundant) global
 constraint.  In contrast to the method described in
 Section~\ref{subsec:confanalysisinflp}, using a dual unbounded ray
 as a set of weights to aggregate  model constraints yields exactly
 one linear constraint.

The proof-constraint along with an activity argument can be used to deduce local lower and upper variable bounds~\cite{achterberg2007constraint}.
Therefore, consider a subproblem with local bounds $[\ell^\prime,u^\prime]$.
For any $i \in \allvars$ with a non-zero coefficient in the proof-constraint the \emph{maximal activity residual} is given by
\begin{align*}
  \maxact{\gamma^tA}{\ell^{\prime}}{u^{\prime}}{i} := \sum_{j \in \allvars\setminus i\colon \gamma^tA_{j} > 0} (\gamma^tA_{j}) u^{\prime}_{j} 
  + \sum_{j \in \allvars\setminus i\colon \gamma^tA_{j} < 0} (\gamma^tA_{j}) \ell^{\prime}_{j},
\end{align*}
\ie the maximal activity over all variables but $x_{i}$.
Hence, valid local bounds are given by
\begin{align*}
  \frac{\gamma^t b - \maxact{\gamma^tA}{\ell^{\prime}}{u^{\prime}}{i}}{a_{i}} 
  \left\{\begin{array}{c} \leq \\ \geq \end{array}\right\}
  x_{i}
  \left\{\begin{array}{c} \text{if } a_i > 0 \\ \text{if } a_i < 0\end{array}\right..
\end{align*}
This is the so-called bound tightening procedure~\cite{brearley75}
which is  widely used in all major MIP solvers, for all kinds of
linear constraints.

Just like the dual ray might be heuristically shrinked to get a short
initial reason for conflict analysis, it might be worthwhile to alter
the proof-constraint itself before using it for propagation. This can
include the application of presolving steps such as coefficient
tightening to the constraint, projecting out continuous variables or
applying mixed-integer rounding to get an alternative globally valid
constraint which might be more powerful to propagate.

Finally, instead of generating a valid constraint from the dual ray, one could
equivalently use the ray itself to simply check for
infeasibility~\cite{polik15ics,polik15ismp} or to estimate the
objective change during branch-and-bound and to derive branching
decisions therefrom. While in Section~\ref{subsec:confanalysisinflp},
we described a way to reduce LP infeasibility analysis to conflict
analysis based on domain propagation, one could as well try to
generate a dual ray by solving the LP relaxation after having
detected infeasibility by propagation.

\section{Managing of Conflicts in a MIP Solver}
\label{sec:conflictpool}

Maintaining and propagating large numbers of conflict constraints
might slow down a solver and create a big burden memory-wise. For
instances with a high throughput of \bandb nodes, a solver like \scip
might easily create hundreds of thousands of conflicts within an hour of running
time. In order to avoid a slowdown or memory short-coming, an aging
mechanism is used within \scip. Once again, aging is a concept inspired
by SAT and CP solving. Every time a conflict constraint is considered
for domain propagation an age counter (individually for each constraint) is increased if no deduction was found.
If a deduction is found, the age will be reset to 0.
If the age reaches a predefined threshold the conflict constraint is permanently deleted.

In SAT and CP, this mechanism is a well-established method to drop conflict constraints
that do not frequently propagate.
In the case of MIP solving, there are two main differences concerning
the \bandb search. First, domain propagation is most often not the
most expensive part of node processing. Second, SAT and CP
solvers often use a pure depth-first-search (DFS) node selection, while
state-of-the-art MIP solvers use some hybrid between DFS and
best-estimate-search or best-first-search (see,
\eg~\cite{achterberg2007constraint,benichou1971experiments,linderoth1999computational}). 
Therefore, it frequently happens that the node processed next is picked from a different part of the tree.

In the following, we describe a pool-based approach to manage conflict constraints.
Here, a pool refers to a fixed-size array that allows direct access to
a particular element and which is independent of the model itself.
The \emph{conflict pool} is used to manage all conflict constraints,
independently whether they were  derived from domain propagation or
an infeasible LP relaxation.
The number of constraints that can be stored within the conflict pool at the same is limited.
In our implementation the maximal size of the conflict pool depends on the number of variables and constraints of the presolved problem.
However, the pool provides space for at least $1\,000$ and at most $50\,000$ conflict constraints at the same time.
The conflict pool allows a central management of conflict constraints independently from the model constraints,
\ie they can be propagated, checked or deleted separately, without
the need to traverse through all constraints.

To drop conflict constraints that don't lead to deductions frequently
we implemented an update-routine that checks the conflict pool
regularly, \eg any time we create the first new conflict at a node.
Moreover, we still use the concept of aging to determine the conflict
constraints that are rarely used in propagation.  Within this update
procedure the oldest conflict constraints are removed.

Beside of the regular checks, the conflict pool is updated every time a new improving incumbent solution is found.
Conflict constraints might depend on a (previous) best known solution,
\eg when the conflict was created from an LP whose infeasibility proof
contained the objective cutoff. Such conflicts become weaker whenever a
new incumbent is found and the chance that they lead to deductions becomes
smaller the more the incumbent improves.
Due to this, for each conflict constraint involving an incumbent solution we store the corresponding objective value.
If this value is sufficiently worse than the new objective value, the
conflict constraint will be permanently deleted.
In our computational experiments (cf. Section~\ref{sec:experiments}) we use a threshold of $5\%$.

\section{Computational Experiments}
\label{sec:experiments}

In our computational experiments,
we compare combinations of the techniques presented in this paper: conflict analysis
and dual ray analysis. To the best of our knowledge, most major MIP
solvers either use conflict analysis of infeasible LPs and domain propagation
(\eg \scip, \textsf{FICO Xpress-Optimizer}) or they employ
 dual ray analysis (\eg \textsf{Gurobi}, \textsf{SAS}).
We will refer to the former as the \conflict setting and to
the latter as the \dualray setting. We compare those to a
setting that uses conflict analysis and dual ray analysis
simultaneously, the \combined setting. Finally, we consider an
extension of the \combined setting that uses a pool for
conflict management, the setting \combinedpool.

All experiments were performed with the non-commercial MIP solver
\scip~\cite{GamrathFischerGallyetal2016} (git hash 60f49ab, based on
\scip~3.2.1.2), using \soplex 2.2.1.3 as LP solver. The
experiments were run on a cluster of identical machines,
each with an Intel Xeon Quad-Core with 3.2\,GHz and 48\,GB of RAM; a
time limit of $3600$ seconds was set.\\

We used two test sets: the \Miplib2010~\cite{KochEtAl2011} benchmark test set and 
a selection of instances taken from the \textsc{Miplib}~\cite{BixbyBoydIndovina1992}, \textsc{Miplib2003}~\cite{AchterbergKochMartin2006},  \Miplib2010,
 the \textsc{Cor@l}~\cite{linderoth2005noncommercial}  collection, the \textsc{Alu}\footnote{The instances are part of the contributed section of \textsc{Miplib2003}},
and the \textsc{Markshare}~\cite{cornuejols1998class} test sets.
From these we selected all instances for which (i) all of the above
settings need at least $100$ nodes, (ii) at least one setting finishes
within the time limit of $3600$ seconds, and (iii) at least one setting analyzes more than $100$ infeasible subproblems successfully.
We refer to this test set as the \Conflict set, since it was
designed to contain instances for which conflict or dual ray analysis is
frequently used.

Aggregated results on the number of generated nodes and needed solving
time can be found in Table~\ref{tab:miplib2010_conflict_aggregated}.
\ifthenelse{\zibreport = 1}%
{%
Detailed results can be found in Table~\ref{tab:miplib_default_dualray_combined_pool_full} and~\ref{tab:conflict_default_dualray_onlyprop_onlyinflp_full} in the appendix.\\
}%
{
Detailed results for each instance and test set separately can be found in the appendix of~\cite{WitzigBertholdHeinz2016}.
}%
We use the \conflict setting as a base line (since it used to
be the \scip default), for which we give actual means of \bandb nodes
and the solving time. For all other settings, we instead give factors
w.r.t.\ the base line. A number greater then one implies that the
setting is inferior and a number less than one implies that the
setting is superior to the \conflict setting.

First of all, we observe that solely using dual ray analysis is
inferior to using conflict analysis on both test sets and w.r.t.\ both
performance measures. Note that we used a basic implementation of dual
ray analysis; a solver that solely relies on it might implement
further extensions that decrease this difference in performance, see
also Section~\ref{sec:outlook}.
However, the combination of conflict and dual ray analysis showed some
significant performance improvements. We observed a speed-up of $3\%$ and $18\%$ on \Miplib2010 and \Conflict, respectively.
Moreover, the number of generated nodes could be reduced by $5\%$ and
$25\%$, respectively. Finally, on the \Conflict test set, the
\combined setting solved one instance more than the
\conflict setting and five more than the \dualray
setting. We take those results as an indicator that the
two techniques complement each other nicely.
In an additional experiment, we also tested to apply conflict analysis
solely from domain propagation or solely from infeasible LPs. Both
variants were inferior to the \conflict setting and are
therefore not discussed in detail%
\ifthenelse{\zibreport = 1}%
{%
\xspace(cf.\ Table~\ref{tab:miplib_default_dualray_onlyprop_onlyinflp_full} 
and~\ref{tab:conflict_default_dualray_onlyprop_onlyinflp_full} in the appendix).\\
}%
{. Detailed results for each instance and test set separately can be found in the appendix of~\cite{WitzigBertholdHeinz2016}.}%

To partially explain the different extent of the improvements on both
tests set, we would like to point out that in the \Miplib2010
benchmark set, there are only 31 instances which fulfill the filtering
criteria mentioned above for the \Conflict set. On those, the
\combined setting is $7.2\%$ faster and needs $15.6\%$ less
nodes than the \conflict setting.

Looking at individual instances, there are a few cases for which the
\combined setting is the clear winner, \eg
\texttt{neos-849702} or\texttt{bnatt350}.
For \texttt{neos-849702} and \texttt{bnatt350}, the \dualray
setting has a timeout, while the \conflict setting is a factor
of  $6.2$ and $1.83$ slower, respectively, than the \combined
setting. At the same time, \texttt{ns1766074} shows the largest
deterioration from using a \combined setting, being a factor
of $1.63$ slower than \conflict and a factor of $1.06$ slower
than \dualray.

As can be seen in Table~\ref{tab:miplib2010_conflict_aggregated},
using a conflict pool in addition to an aging system makes hardly any
difference w.r.t.\ the overall performance.

\begin{table}[t]
\begin{centering}
\small
\setlength{\tabcolsep}{2.5pt}
\begin{tabularx}{\textwidth}{lrrrrrrrrrrrr}
\toprule
  & \multicolumn{3}{c}{\conflict} & \multicolumn{3}{c}{\dualray} & \multicolumn{3}{c}{\combined} & \multicolumn{3}{c}{\combinedpool} \\ 
 \cmidrule(lr){2-4} \cmidrule(lr){5-7} \cmidrule(lr){8-10} \cmidrule(lr){11-13}
Test set & \multicolumn{1}{c}{$\#$} & \multicolumn{1}{c}{$n$} & \multicolumn{1}{c}{$t$}  & \multicolumn{1}{c}{$\#$} & \multicolumn{1}{c}{$n_{Q}$} & \multicolumn{1}{c}{$t_{Q}$}  & \multicolumn{1}{c}{$\#$} & \multicolumn{1}{c}{$n_{Q}$} & \multicolumn{1}{c}{$t_{Q}$}  & \multicolumn{1}{c}{$\#$} & \multicolumn{1}{c}{$n_{Q}$} & \multicolumn{1}{c}{$t_{Q}$}  \\ 
\midrule
\textsc{Miplib2010}  &   60 &     14382 &       686 &   57 & \textit{\textcolor{ red}{     1.365}} & \textit{\textcolor{ red}{     1.167}} &   60 &                               0.955   &                               0.977   &   60 &                               0.957  &                               0.975   \\
\textsc{Conflict}  &  105 &     16769 &       143 &  101 & \textit{\textcolor{ red}{     1.616}} & \textit{\textcolor{ red}{     1.256}} &  106 & \textbf{\textcolor{blue}{     0.755}} & \textbf{\textcolor{blue}{     0.827}} &  106 & \textbf{\textcolor{blue}{     0.759}} & \textbf{\textcolor{blue}{     0.829}} \\
\bottomrule
\end{tabularx}
\caption{Aggregated computational results. Columns marked with $\#$ show the number of solved instances. Columns 3 and 4 show the shifted geometric mean of absolute numbers of
          generated nodes ($n$, shift = $100$) and needed solving time in seconds ($t$, shift = $10$), respectively.
          All remaining columns show the relative number of generated nodes ($n_{Q}$) and needed solving time ($t_{Q}$) w.r.t. Column~3 and~4, respectively.}
\label{tab:miplib2010_conflict_aggregated}
\end{centering}
\end{table}

\section{Conclusion and Outlook}
\label{sec:outlook}

In this short paper we discussed the similarities and differences of conflict analysis and dual ray analysis in solving MIPs.
Our computational study indicates that a combination of both
approaches can enhance the performance of a state-of-the-art MIP
solver significantly.
On instances where the analysis of infeasible subproblems succeeds frequently, the solving
time improved by $17.3\%$ and the number of \bandb nodes by $24.5\%$.
In contrast to that, using a pool-based approach in addition to an
aging mechanism to manage conflict constraints showed hardly any impact.\\

There are several instances for which using either dual ray analysis
or conflict analysis exclusively outperformed the combination of both.
Thus, we will plan to investigate a dynamic mechanism to switch between both techniques.
Furthermore, applying dual ray analysis for infeasibility deduced by
domain propagation as well as using more preprocessing (\eg mixed integer
rounding, projecting out continuous variables, etc.) techniques to
modify constraints derived from dual ray analysis appear as promising
directions for future research.

\subsection*{Acknowledgments}

The work for this article has been conducted within the Research Campus Modal
funded by the German Federal Ministry of Education and Research (fund number 05M14ZAM).

\bibliographystyle{abbrv}
\bibliography{Bibliography}

\begin{thebibliography}{10}

\bibitem{achterberg2007conflict}
T.~Achterberg.
\newblock Conflict analysis in mixed integer programming.
\newblock {\em Discrete Optimization}, 4(1):4--20, 2007.

\bibitem{achterberg2007constraint}
T.~Achterberg.
\newblock Constraint integer programming, 2007.

\bibitem{AchterbergBerthold2009}
T.~Achterberg and T.~Berthold.
\newblock Hybrid branching.
\newblock In W.-J. van Hoeve and J.~N. Hooker, editors, {\em Integration of AI
  and OR Techniques in Constraint Programming for Combinatorial Optimization
  Problems, 6th International Conference, CPAIOR 2009}, volume 5547 of {\em
  Lecture Notes in Computer Science}, pages 309--311. Springer Berlin
  Heidelberg, May 2009.

\bibitem{AchterbergKochMartin2006}
T.~Achterberg, T.~Koch, and A.~Martin.
\newblock {MIPLIB} 2003.
\newblock {\em Operations Research Letters}, 34(4):361--372, 2006.

\bibitem{benichou1971experiments}
M.~B{\'e}nichou, J.-M. Gauthier, P.~Girodet, G.~Hentges, G.~Ribi{\`e}re, and
  O.~Vincent.
\newblock Experiments in mixed-integer linear programming.
\newblock {\em Mathematical Programming}, 1(1):76--94, 1971.

\bibitem{BertholdGleixner2013}
T.~Berthold and A.~M. Gleixner.
\newblock Undercover: a primal {MINLP} heuristic exploring a largest sub-{MIP}.
\newblock {\em Mathematical Programming}, 144(1--2):315--346, 2014.

\bibitem{BertholdHendel2015}
T.~Berthold and G.~Hendel.
\newblock Shift-and-propagate.
\newblock {\em Journal of Heuristics}, 21(1):73--106, 2015.

\bibitem{BixbyBoydIndovina1992}
R.~E. Bixby, E.~A. Boyd, and R.~R. Indovina.
\newblock {MIPLIB}: {A} test set of mixed integer programming problems.
\newblock {\em SIAM News}, 25:16, 1992.

\bibitem{brearley75}
A.~Brearley, G.~Mitra, and H.~Williams.
\newblock Analysis of mathematical programming problems prior to applying the
  simplex algorithm.
\newblock {\em Mathematical Programming}, 8:54--83, 1975.

\bibitem{cornuejols1998class}
G.~Cornu{\'e}jols and M.~Dawande.
\newblock A class of hard small 0-1 programs.
\newblock In {\em International Conference on Integer Programming and
  Combinatorial Optimization}, pages 284--293. Springer, 1998.

\bibitem{Dakin1965}
R.~J. Dakin.
\newblock A tree-search algorithm for mixed integer programming problems.
\newblock {\em The Computer Journal}, 8(3):250--255, 1965.

\bibitem{DaveyBolandStuckey2002}
B.~Davey, N.~Boland, and P.~J. Stuckey.
\newblock Efficient intelligent backtracking using linear programming.
\newblock {\em INFORMS Journal of Computing}, 14(4):373--386, 2002.

\bibitem{GamrathFischerGallyetal2016}
G.~Gamrath, T.~Fischer, T.~Gally, A.~M. Gleixner, G.~Hendel, T.~Koch, S.~J.
  Maher, M.~Miltenberger, B.~M{\"u}ller, M.~E. Pfetsch, C.~Puchert,
  D.~Rehfeldt, S.~Schenker, R.~Schwarz, F.~Serrano, Y.~Shinano, S.~Vigerske,
  D.~Weninger, M.~Winkler, J.~T. Witt, and J.~Witzig.
\newblock The scip optimization suite 3.2.
\newblock Technical Report 15-60, ZIB, Takustr.7, 14195 Berlin, 2016.

\bibitem{ginsberg1993dynamic}
M.~L. Ginsberg.
\newblock Dynamic backtracking.
\newblock {\em Journal of Artificial Intelligence Research}, pages 25--46,
  1993.

\bibitem{jiang1994no}
Y.~Jiang, T.~Richards, and B.~Richards.
\newblock No-good backmarking with min-con ict repair in constraint
  satisfaction and optimization.
\newblock In {\em PPCP}, volume~94, pages 2--4. Citeseer, 1994.

\bibitem{KilincNemhauserSavelsbergh2010}
F.~K{\i}l{\i}n\c{c}~Karzan, G.~L. Nemhauser, and M.~W.~P. Savelsbergh.
\newblock Information-based branching schemes for binary linear mixed-integer
  programs.
\newblock {\em Mathematical Programming Computation}, 1(4):249--293, 2009.

\bibitem{KochEtAl2011}
T.~Koch, T.~Achterberg, E.~Andersen, O.~Bastert, T.~Berthold, R.~E. Bixby,
  E.~Danna, G.~Gamrath, A.~M. Gleixner, S.~Heinz, A.~Lodi, H.~Mittelmann,
  T.~Ralphs, D.~Salvagnin, D.~E. Steffy, and K.~Wolter.
\newblock {MIPLIB} 2010.
\newblock {\em Mathematical Programming Computation}, 3(2):103--163, 2011.

\bibitem{LandDoig1960}
A.~H. Land and A.~G. Doig.
\newblock An automatic method of solving discrete programming problems.
\newblock {\em Econometrica}, 28(3):497--520, 1960.

\bibitem{linderoth2005noncommercial}
J.~T. Linderoth and T.~K. Ralphs.
\newblock Noncommercial software for mixed-integer linear programming.
\newblock {\em Integer programming: theory and practice}, 3:253--303, 2005.

\bibitem{linderoth1999computational}
J.~T. Linderoth and M.~W. Savelsbergh.
\newblock A computational study of search strategies for mixed integer
  programming.
\newblock {\em INFORMS Journal on Computing}, 11(2):173--187, 1999.

\bibitem{marques1999grasp}
J.~P. Marques-Silva and K.~Sakallah.
\newblock Grasp: A search algorithm for propositional satisfiability.
\newblock {\em Computers, IEEE Transactions on}, 48(5):506--521, 1999.

\bibitem{polik15ics}
I.~P\'olik.
\newblock (re)using dual information in milp.
\newblock In {\em INFORMS Computing Society conference}, Richmond, VA, 2015.

\bibitem{polik15ismp}
I.~P\'olik.
\newblock Some more ways to use dual information in milp.
\newblock In {\em International Symposium on Mathematical Programming},
  Pittsburgh, PA, 2015.

\bibitem{SandholmShields2006}
T.~Sandholm and R.~Shields.
\newblock Nogood learning for mixed integer programming.
\newblock In {\em Workshop on Hybrid Methods and Branching Rules in
  Combinatorial Optimization, Montr{\'e}al}, 2006.

\bibitem{stallman1977forward}
R.~M. Stallman and G.~J. Sussman.
\newblock Forward reasoning and dependency-directed backtracking in a system
  for computer-aided circuit analysis.
\newblock {\em Artificial intelligence}, 9(2):135--196, 1977.

\end{thebibliography}

\ifthenelse{\zibreport=1}
{
\clearpage
\pagebreak

\begin{appendix}

\section{Appendix}

A detailed overview of all computational results on \Miplib2010 and \Conflict test set can be found in 
Table~\ref{tab:miplib_default_dualray_combined_pool_full} --~\ref{tab:conflict_default_dualray_onlyprop_onlyinflp_full}.
For each table we use the \conflict setting as a base line, for which we give actual means of \bandb nodes and the solving time.
For all other settings, we instead give factors w.r.t.\ the base line.
A number greater then one implies that the setting is inferior and a number less than one implies that the setting is superior to the \conflict setting.

A comparison between the \conflict, \dualray, \combined, and \combinedpool setting can be found in Table~\ref{tab:miplib_default_dualray_combined_pool_full}
and~\ref{tab:conflict_default_dualray_combined_pool_full}.

In addition, results for applying conflict analysis solely from domain propagation or solely from infeasible LPs can be found in 
Table~\ref{tab:miplib_default_dualray_onlyprop_onlyinflp_full} and~\ref{tab:conflict_default_dualray_onlyprop_onlyinflp_full}.

\begin{centering}
\footnotesize
\begin{landscape}
\setlength{\tabcolsep}{3pt}
\begin{longtable}{lrrrrrrrrrrrrrrrr}
\caption{Detailed computational results on \Miplib2010 test set. The table shows the absolute number of generated nodes ($n$) and needed solving time in seconds ($t$),
          as well as the relative number of generated nodes ($n_{Q}$) and needed solving time ($t_{Q}$) w.r.t. Columns~2 and~3 (\conflict).
          All changes in the number of nodes or solving time of at least $5\%$ are highlighted in bold and blue (\textbf{\textcolor{blue}{improvement}}) and
          italic and red (\textit{\textcolor{red}{deterioration}}).}\label{tab:miplib_default_dualray_combined_pool_full}\\ 
\toprule
  & \multicolumn{2}{c}{\conflict} & \multicolumn{4}{c}{\dualray} & \multicolumn{4}{c}{\combined} & \multicolumn{4}{c}{\combinedpool} \\ 
 \cmidrule(lr){2-3} \cmidrule(lr){4-7} \cmidrule(lr){8-11} \cmidrule(lr){12-15}
Instance & \multicolumn{1}{c}{$n$} & \multicolumn{1}{c}{$t$}  & \multicolumn{1}{c}{$n$} & \multicolumn{1}{c}{$n_{Q}$} & \multicolumn{1}{c}{$t$} & \multicolumn{1}{c}{$t_{Q}$}  & \multicolumn{1}{c}{$n$} & \multicolumn{1}{c}{$n_{Q}$} & \multicolumn{1}{c}{$t$} & \multicolumn{1}{c}{$t_{Q}$}  & \multicolumn{1}{c}{$n$} & \multicolumn{1}{c}{$n_{Q}$} & \multicolumn{1}{c}{$t$} & \multicolumn{1}{c}{$t_{Q}$}  \\ 
\midrule
\endfirsthead

\caption{Detailed computational results on \Miplib2010 test set. The table shows the absolute number of generated nodes ($n$) and needed solving time in seconds ($t$),
          as well as the relative number of generated nodes ($n_{Q}$) and needed solving time ($t_{Q}$) w.r.t. Columns~2 and~3 (\conflict).
          All changes in the number of nodes or solving time of at least $5\%$ are highlighted in bold and blue (\textbf{\textcolor{blue}{improvement}}) and
          italic and red (\textit{\textcolor{red}{deterioration}}).}\\ 
\toprule
  & \multicolumn{2}{c}{\conflict} & \multicolumn{4}{c}{\dualray} & \multicolumn{4}{c}{\combined} & \multicolumn{4}{c}{\combinedpool} \\ 
 \cmidrule(lr){2-3} \cmidrule(lr){4-7} \cmidrule(lr){8-11} \cmidrule(lr){12-15}
Instance & \multicolumn{1}{c}{$n$} & \multicolumn{1}{c}{$t$}  & \multicolumn{1}{c}{$n$} & \multicolumn{1}{c}{$n_{Q}$} & \multicolumn{1}{c}{$t$} & \multicolumn{1}{c}{$t_{Q}$}  & \multicolumn{1}{c}{$n$} & \multicolumn{1}{c}{$n_{Q}$} & \multicolumn{1}{c}{$t$} & \multicolumn{1}{c}{$t_{Q}$}  & \multicolumn{1}{c}{$n$} & \multicolumn{1}{c}{$n_{Q}$} & \multicolumn{1}{c}{$t$} & \multicolumn{1}{c}{$t_{Q}$}  \\ 
\midrule
\endhead

\midrule
\multicolumn{15}{r}{cont. on next page}\\
\bottomrule
\endfoot

\bottomrule
\endlastfoot
30n20b8                   &        51 &       268 &        933 & \textit{\textcolor{ red}{    18.294}} &    1591.27 & \textit{\textcolor{ red}{     5.933}} &         51 &                               1.000   &     269.69 &                               1.006   &         51 &                               1.000   &     269.43 &                               1.005   \\
acc-tight5                &       920 &       140 &        357 & \textbf{\textcolor{blue}{     0.388}} &      76.88 & \textbf{\textcolor{blue}{     0.551}} &       1107 & \textit{\textcolor{ red}{     1.203}} &     153.85 & \textit{\textcolor{ red}{     1.102}} &       1107 & \textit{\textcolor{ red}{     1.203}} &     153.74 & \textit{\textcolor{ red}{     1.101}} \\
aflow40b                  &    163981 &      1333 &     160367 &                               0.978   &     888.21 & \textbf{\textcolor{blue}{     0.666}} &     135529 & \textbf{\textcolor{blue}{     0.826}} &     634.29 & \textbf{\textcolor{blue}{     0.476}} &     135018 & \textbf{\textcolor{blue}{     0.823}} &     631.47 & \textbf{\textcolor{blue}{     0.474}} \\
air04                     &       218 &        74 &        180 & \textbf{\textcolor{blue}{     0.826}} &      73.71 &                               0.997   &        218 &                               1.000   &      69.62 & \textbf{\textcolor{blue}{     0.942}} &        218 &                               1.000   &      70.40 &                               0.952   \\
app1-2                    &         2 &      3600 &         71 & \textit{\textcolor{ red}{    35.500}} &    2100.69 & \textbf{\textcolor{blue}{     0.584}} &          2 &                               1.000   &    3600.00 &                               1.000   &          2 &                               1.000   &    3600.00 &                               1.000   \\
\end{longtable}
\end{landscape}
\end{centering}

\begin{centering}
\footnotesize
\begin{landscape}
\setlength{\tabcolsep}{3pt}
\begin{longtable}{lrrrrrrrrrrrrrrrr}
\caption{Detailed computational results on \Conflict test set. The table shows the absolute number of generated nodes ($n$) and needed solving time in seconds ($t$),
          as well as the relative number of generated nodes ($n_{Q}$) and needed solving time ($t_{Q}$) w.r.t. Columns~2 and~3 (\conflict).
          All changes in the number of nodes or solving time of at least $5\%$ are highlighted in bold and blue (\textbf{\textcolor{blue}{improvement}}) and
          italic and red (\textit{\textcolor{red}{deterioration}}).}\label{tab:conflict_default_dualray_combined_pool_full}\\ 
\toprule
  & \multicolumn{2}{c}{\conflict} & \multicolumn{4}{c}{\dualray} & \multicolumn{4}{c}{\combined} & \multicolumn{4}{c}{\combinedpool} \\ 
 \cmidrule(lr){2-3} \cmidrule(lr){4-7} \cmidrule(lr){8-11} \cmidrule(lr){12-15}
Instance & \multicolumn{1}{c}{$n$} & \multicolumn{1}{c}{$t$}  & \multicolumn{1}{c}{$n$} & \multicolumn{1}{c}{$n_{Q}$} & \multicolumn{1}{c}{$t$} & \multicolumn{1}{c}{$t_{Q}$}  & \multicolumn{1}{c}{$n$} & \multicolumn{1}{c}{$n_{Q}$} & \multicolumn{1}{c}{$t$} & \multicolumn{1}{c}{$t_{Q}$}  & \multicolumn{1}{c}{$n$} & \multicolumn{1}{c}{$n_{Q}$} & \multicolumn{1}{c}{$t$} & \multicolumn{1}{c}{$t_{Q}$}  \\ 
\midrule
\endfirsthead

\caption{Detailed computational results on \Conflict test set. The table shows the absolute number of generated nodes ($n$) and needed solving time in seconds ($t$),
          as well as the relative number of generated nodes ($n_{Q}$) and needed solving time ($t_{Q}$) w.r.t. Columns~2 and~3 (\conflict).
          All changes in the number of nodes or solving time of at least $5\%$ are highlighted in bold and blue (\textbf{\textcolor{blue}{improvement}}) and
          italic and red (\textit{\textcolor{red}{deterioration}}).}\\ 
\toprule
  & \multicolumn{2}{c}{\conflict} & \multicolumn{4}{c}{\dualray} & \multicolumn{4}{c}{\combined} & \multicolumn{4}{c}{\combinedpool} \\ 
 \cmidrule(lr){2-3} \cmidrule(lr){4-7} \cmidrule(lr){8-11} \cmidrule(lr){12-15}
Instance & \multicolumn{1}{c}{$n$} & \multicolumn{1}{c}{$t$}  & \multicolumn{1}{c}{$n$} & \multicolumn{1}{c}{$n_{Q}$} & \multicolumn{1}{c}{$t$} & \multicolumn{1}{c}{$t_{Q}$}  & \multicolumn{1}{c}{$n$} & \multicolumn{1}{c}{$n_{Q}$} & \multicolumn{1}{c}{$t$} & \multicolumn{1}{c}{$t_{Q}$}  & \multicolumn{1}{c}{$n$} & \multicolumn{1}{c}{$n_{Q}$} & \multicolumn{1}{c}{$t$} & \multicolumn{1}{c}{$t_{Q}$}  \\ 
\midrule
\endhead

\midrule
\multicolumn{15}{r}{cont. on next page}\\
\bottomrule
\endfoot

\bottomrule
\endlastfoot
acc-tight5                &       920 &       139 &        357 & \textbf{\textcolor{blue}{     0.388}} &      80.27 & \textbf{\textcolor{blue}{     0.577}} &       1107 & \textit{\textcolor{ red}{     1.203}} &     153.52 & \textit{\textcolor{ red}{     1.104}} &       1107 & \textit{\textcolor{ red}{     1.203}} &     154.34 & \textit{\textcolor{ red}{     1.109}} \\
aflow40b                  &    163981 &      1338 &     160367 &                               0.978   &     945.95 & \textbf{\textcolor{blue}{     0.707}} &     135529 & \textbf{\textcolor{blue}{     0.826}} &     649.99 & \textbf{\textcolor{blue}{     0.486}} &     135018 & \textbf{\textcolor{blue}{     0.823}} &     632.94 & \textbf{\textcolor{blue}{     0.473}} \\
alu4\_7                   &       279 &         1 &       1643 & \textit{\textcolor{ red}{     5.889}} &       1.88 & \textit{\textcolor{ red}{     1.709}} &        250 & \textbf{\textcolor{blue}{     0.896}} &       1.12 &                               1.018   &        250 & \textbf{\textcolor{blue}{     0.896}} &       1.17 & \textit{\textcolor{ red}{     1.064}} \\
alu4\_8                   &       649 &         1 &      28296 & \textit{\textcolor{ red}{    43.599}} &      10.73 & \textit{\textcolor{ red}{     8.007}} &        710 & \textit{\textcolor{ red}{     1.094}} &       1.41 & \textit{\textcolor{ red}{     1.052}} &        710 & \textit{\textcolor{ red}{     1.094}} &       1.41 & \textit{\textcolor{ red}{     1.052}} \\
alu5\_7                   &      1704 &         2 &       6467 & \textit{\textcolor{ red}{     3.795}} &       3.97 & \textit{\textcolor{ red}{     2.005}} &       1368 & \textbf{\textcolor{blue}{     0.803}} &       1.67 & \textbf{\textcolor{blue}{     0.843}} &       1368 & \textbf{\textcolor{blue}{     0.803}} &       1.66 & \textbf{\textcolor{blue}{     0.838}} \\
\end{longtable}
\end{landscape}
\end{centering}

\begin{centering}
\footnotesize
\begin{landscape}
\setlength{\tabcolsep}{2.5pt}
\begin{longtable}{lrrrrrrrrrrrrrrrr}
\caption{Detailed computational results on \Miplib2010 test set. The table shows the absolute number of generated nodes ($n$) and needed solving time in seconds ($t$),
          as well as the relative number of generated nodes ($n_{Q}$) and needed solving time ($t_{Q}$) w.r.t. Columns~2 and~3 (\conflict).
          All changes in the number of nodes or solving time of at least $5\%$ are highlighted in bold and blue (\textbf{\textcolor{blue}{improvement}}) and
          italic and red (\textit{\textcolor{red}{deterioration}}).}\label{tab:miplib_default_dualray_onlyprop_onlyinflp_full}\\ 
\toprule
  & \multicolumn{2}{c}{\conflict} & \multicolumn{4}{c}{\dualray} & \multicolumn{4}{c}{\texttt{inf. LP}} & \multicolumn{4}{c}{\texttt{dom. prop.}} \\ 
 \cmidrule(lr){2-3} \cmidrule(lr){4-7} \cmidrule(lr){8-11} \cmidrule(lr){12-15}
Instance & \multicolumn{1}{c}{$n$} & \multicolumn{1}{c}{$t$}  & \multicolumn{1}{c}{$n$} & \multicolumn{1}{c}{$n_{Q}$} & \multicolumn{1}{c}{$t$} & \multicolumn{1}{c}{$t_{Q}$}  & \multicolumn{1}{c}{$n$} & \multicolumn{1}{c}{$n_{Q}$} & \multicolumn{1}{c}{$t$} & \multicolumn{1}{c}{$t_{Q}$}  & \multicolumn{1}{c}{$n$} & \multicolumn{1}{c}{$n_{Q}$} & \multicolumn{1}{c}{$t$} & \multicolumn{1}{c}{$t_{Q}$}  \\ 
\midrule
\endfirsthead

\caption{Detailed computational results on \Miplib2010 test set. The table shows the absolute number of generated nodes ($n$) and needed solving time in seconds ($t$),
          as well as the relative number of generated nodes ($n_{Q}$) and needed solving time ($t_{Q}$) w.r.t. Columns~2 and~3 (\conflict).
          All changes in the number of nodes or solving time of at least $5\%$ are highlighted in bold and blue (\textbf{\textcolor{blue}{improvement}}) and
          italic and red (\textit{\textcolor{red}{deterioration}}).}\\ 
\toprule
  & \multicolumn{2}{c}{\conflict} & \multicolumn{4}{c}{\dualray} & \multicolumn{4}{c}{\texttt{inf. LP}} & \multicolumn{4}{c}{\texttt{dom. prop.}} \\ 
 \cmidrule(lr){2-3} \cmidrule(lr){4-7} \cmidrule(lr){8-11} \cmidrule(lr){12-15}
Instance & \multicolumn{1}{c}{$n$} & \multicolumn{1}{c}{$t$}  & \multicolumn{1}{c}{$n$} & \multicolumn{1}{c}{$n_{Q}$} & \multicolumn{1}{c}{$t$} & \multicolumn{1}{c}{$t_{Q}$}  & \multicolumn{1}{c}{$n$} & \multicolumn{1}{c}{$n_{Q}$} & \multicolumn{1}{c}{$t$} & \multicolumn{1}{c}{$t_{Q}$}  & \multicolumn{1}{c}{$n$} & \multicolumn{1}{c}{$n_{Q}$} & \multicolumn{1}{c}{$t$} & \multicolumn{1}{c}{$t_{Q}$}  \\ 
\midrule
\endhead

\midrule
\multicolumn{15}{r}{cont. on next page}\\
\bottomrule
\endfoot

\bottomrule
\endlastfoot
30n20b8                   &        51 &       268 &        933 & \textit{\textcolor{ red}{    18.294}} &    1591.27 & \textit{\textcolor{ red}{     5.933}} &       1207 & \textit{\textcolor{ red}{    23.667}} &    1477.92 & \textit{\textcolor{ red}{     5.511}} &       6754 & \textit{\textcolor{ red}{   132.431}} &    3600.00 & \textit{\textcolor{ red}{    13.423}} \\
acc-tight5                &       920 &       140 &        357 & \textbf{\textcolor{blue}{     0.388}} &      76.88 & \textbf{\textcolor{blue}{     0.551}} &        207 & \textbf{\textcolor{blue}{     0.225}} &      69.40 & \textbf{\textcolor{blue}{     0.497}} &       1615 & \textit{\textcolor{ red}{     1.755}} &     280.53 & \textit{\textcolor{ red}{     2.009}} \\
aflow40b                  &    163981 &      1333 &     160367 &                               0.978   &     888.21 & \textbf{\textcolor{blue}{     0.666}} &     158904 &                               0.969   &    1259.54 & \textbf{\textcolor{blue}{     0.945}} &     199108 & \textit{\textcolor{ red}{     1.214}} &    1486.33 & \textit{\textcolor{ red}{     1.115}} \\
air04                     &       218 &        74 &        180 & \textbf{\textcolor{blue}{     0.826}} &      73.71 &                               0.997   &        180 & \textbf{\textcolor{blue}{     0.826}} &      73.85 &                               0.999   &        218 &                               1.000   &      69.87 & \textbf{\textcolor{blue}{     0.945}} \\
app1-2                    &         2 &      3600 &         71 & \textit{\textcolor{ red}{    35.500}} &    2100.69 & \textbf{\textcolor{blue}{     0.584}} &          2 &                               1.000   &    3600.00 &                               1.000   &          4 & \textit{\textcolor{ red}{     2.000}} &    3600.00 &                               1.000   \\
\end{longtable}
\end{landscape}
\end{centering}

\begin{centering}
\footnotesize
\begin{landscape}
\setlength{\tabcolsep}{2.5pt}
\begin{longtable}{lrrrrrrrrrrrrrrrr}
\caption{Detailed computational results on \Conflict test set. The table shows the absolute number of generated nodes ($n$) and needed solving time in seconds ($t$),
          as well as the relative number of generated nodes ($n_{Q}$) and needed solving time ($t_{Q}$) w.r.t. Columns~2 and~3 (\conflict).
          All changes in the number of nodes or solving time of at least $5\%$ are highlighted in bold and blue (\textbf{\textcolor{blue}{improvement}}) and
          italic and red (\textit{\textcolor{red}{deterioration}}).}\label{tab:conflict_default_dualray_onlyprop_onlyinflp_full}\\ 
\toprule
  & \multicolumn{2}{c}{\conflict} & \multicolumn{4}{c}{\dualray} & \multicolumn{4}{c}{\texttt{inf. LP}} & \multicolumn{4}{c}{\texttt{dom. prop.}} \\ 
 \cmidrule(lr){2-3} \cmidrule(lr){4-7} \cmidrule(lr){8-11} \cmidrule(lr){12-15}
Instance & \multicolumn{1}{c}{$n$} & \multicolumn{1}{c}{$t$}  & \multicolumn{1}{c}{$n$} & \multicolumn{1}{c}{$n_{Q}$} & \multicolumn{1}{c}{$t$} & \multicolumn{1}{c}{$t_{Q}$}  & \multicolumn{1}{c}{$n$} & \multicolumn{1}{c}{$n_{Q}$} & \multicolumn{1}{c}{$t$} & \multicolumn{1}{c}{$t_{Q}$}  & \multicolumn{1}{c}{$n$} & \multicolumn{1}{c}{$n_{Q}$} & \multicolumn{1}{c}{$t$} & \multicolumn{1}{c}{$t_{Q}$}  \\ 
\midrule
\endfirsthead

\caption{Detailed computational results on \Conflict test set. The table shows the absolute number of generated nodes ($n$) and needed solving time in seconds ($t$),
          as well as the relative number of generated nodes ($n_{Q}$) and needed solving time ($t_{Q}$) w.r.t. Columns~2 and~3 (\conflict).
          All changes in the number of nodes or solving time of at least $5\%$ are highlighted in bold and blue (\textbf{\textcolor{blue}{improvement}}) and
          italic and red (\textit{\textcolor{red}{deterioration}}).}\\ 
\toprule
  & \multicolumn{2}{c}{\conflict} & \multicolumn{4}{c}{\dualray} & \multicolumn{4}{c}{\texttt{inf. LP}} & \multicolumn{4}{c}{\texttt{dom. prop.}} \\ 
 \cmidrule(lr){2-3} \cmidrule(lr){4-7} \cmidrule(lr){8-11} \cmidrule(lr){12-15}
Instance & \multicolumn{1}{c}{$n$} & \multicolumn{1}{c}{$t$}  & \multicolumn{1}{c}{$n$} & \multicolumn{1}{c}{$n_{Q}$} & \multicolumn{1}{c}{$t$} & \multicolumn{1}{c}{$t_{Q}$}  & \multicolumn{1}{c}{$n$} & \multicolumn{1}{c}{$n_{Q}$} & \multicolumn{1}{c}{$t$} & \multicolumn{1}{c}{$t_{Q}$}  & \multicolumn{1}{c}{$n$} & \multicolumn{1}{c}{$n_{Q}$} & \multicolumn{1}{c}{$t$} & \multicolumn{1}{c}{$t_{Q}$}  \\ 
\midrule
\endhead

\midrule
\multicolumn{15}{r}{cont. on next page}\\
\bottomrule
\endfoot

\bottomrule
\endlastfoot
acc-tight5                &       920 &       139 &        357 & \textbf{\textcolor{blue}{     0.388}} &      80.27 & \textbf{\textcolor{blue}{     0.577}} &        207 & \textbf{\textcolor{blue}{     0.225}} &      66.28 & \textbf{\textcolor{blue}{     0.476}} &       1615 & \textit{\textcolor{ red}{     1.755}} &     280.88 & \textit{\textcolor{ red}{     2.019}} \\
aflow40b                  &    163981 &      1338 &     160367 &                               0.978   &     945.95 & \textbf{\textcolor{blue}{     0.707}} &     158904 &                               0.969   &    1262.21 & \textbf{\textcolor{blue}{     0.943}} &     199108 & \textit{\textcolor{ red}{     1.214}} &    1486.37 & \textit{\textcolor{ red}{     1.111}} \\
alu4\_7                   &       279 &         1 &       1643 & \textit{\textcolor{ red}{     5.889}} &       1.88 & \textit{\textcolor{ red}{     1.709}} &        343 & \textit{\textcolor{ red}{     1.229}} &       1.16 & \textit{\textcolor{ red}{     1.055}} &        301 & \textit{\textcolor{ red}{     1.079}} &       1.04 & \textbf{\textcolor{blue}{     0.945}} \\
alu4\_8                   &       649 &         1 &      28296 & \textit{\textcolor{ red}{    43.599}} &      10.73 & \textit{\textcolor{ red}{     8.007}} &       1802 & \textit{\textcolor{ red}{     2.777}} &       1.72 & \textit{\textcolor{ red}{     1.284}} &       1141 & \textit{\textcolor{ red}{     1.758}} &       1.77 & \textit{\textcolor{ red}{     1.321}} \\
alu5\_7                   &      1704 &         2 &       6467 & \textit{\textcolor{ red}{     3.795}} &       3.97 & \textit{\textcolor{ red}{     2.005}} &       4051 & \textit{\textcolor{ red}{     2.377}} &       2.68 & \textit{\textcolor{ red}{     1.354}} &       2147 & \textit{\textcolor{ red}{     1.260}} &       2.12 & \textit{\textcolor{ red}{     1.071}} \\
\end{longtable}
\end{landscape}
\end{centering}

\end{appendix}
}
{}

\end{document}